# Spectral identification of networks with generalized diffusive coupling

M. Gulina * A. Mauroy *

* Department of Mathematics and Namur Institute for Complex Systems (naXys), University of Namur, 5000 - Belgium

**Abstract:** Spectral network identification aims at inferring the eigenvalues of the Laplacian matrix of a network from measurement data. This allows to capture global information on the network structure from local measurements at a few number of nodes. In this paper, we consider the spectral network identification problem in the generalized setting of a vector-valued diffusive coupling. The feasibility of this problem is investigated and theoretical results on the properties of the associated generalized eigenvalue problem are obtained. Finally, we propose a numerical method to solve the generalized network identification problem, which relies on dynamic mode decomposition and leverages the above theoretical results.

## 1. INTRODUCTION

In some situations, only local measurements of a few nodes of a network are available, while the whole network topology is unknown. Since the network structure affects the dynamics of coupled units attached to the nodes, it is possible to infer some information on the global structure of the network from these local measurements. In particular, this can be done in the context of spectral network identification introduced in Mauroy and Hendrickx (2017). In this setting, the eigenvalues of the Jacobian matrix characterizing the linearized dynamics are computed from local measurement through the so-called DMD algorithm (see Schmid (2010); Tu et al. (2014)) and are shown to be related to the eigenvalues of the Laplacian matrix of the network. Thanks to spectral graph theory, the obtained Laplacian spectrum can in turn be used to infer global structural properties of the network such as the minimum, mean and maximum node degrees (see *e.g.* Fiedler (1973)). However, the results of Mauroy and Hendrickx (2017) are valid only for units interacting through a scalar diffusive coupling. This restrictive condition is not satisfied in many cases, such as reaction-diffusion networks.

This paper aims to extend the applicability of the spectral network identification framework to the setting of vector-valued diffusive coupling between interacting units. In this context, the feasibility problem appears to be more involved than in the scalar-valued coupling case and theoretical results related to a generalized eigenvalue problem are provided. Furthermore, a method for spectral network identification is developed in this generalized case.

The rest of the paper is organized as follows. In Section 2, the spectral network identification problem is introduced in the generalized setting of vector-valued diffusive coupling. The feasibility of this problem is investigated in Section 3 and a numerical method is proposed in Section 4, based on theoretical results related to a generalized eigenvalue problem. The method is illustrated with several examples in Section 5. Finally, concluding remarks and perspectives are discussed in Section 6.

## 2. SPECTRAL NETWORK IDENTIFICATION

We consider a nonlinear dynamical system over a network of $n$ nodes. The unit attached to node $k$ is described by $m$ states whose evolution is governed by the diffusely coupled nonlinear dynamics

$$\begin{cases} \dot{x}_k = F(x_k) + G(x_k)u_k & \in \mathbb{R}^m \\ u_k = \sum_{j=1}^{n} W_{kj}(y_k - y_j) & \in \mathbb{R}^r \\ y_k = H(x_k), & \in \mathbb{R}^r \end{cases} \quad (1)$$

with the (continuously differentiable) functions $F : \mathbb{R}^m \to \mathbb{R}^m$, $G : \mathbb{R}^m \to \mathbb{R}^{m \times r}$ and $H : \mathbb{R}^m \to \mathbb{R}^r$. We assume that the dynamics is so that the units asymptotically reach a synchronized state $x_1 = \cdots = x_n = x^*$, *i.e.*, $\lim_{t\to\infty} x_k(t) = x^*$. Note that the coupling coefficients $W_{kj}$ are the entries of the adjacency matrix $W \in \mathbb{R}^{n \times n}$ of the network. Moreover, the Laplacian matrix is given by $L = D - W$, where $D = \text{diag}(d_1, \ldots, d_n)$ is the degree matrix with $d_i = \sum_{k=1}^{n} W_{ik}$. For the sake of brevity, we will also denote $\mathcal{I}_n = \{1, \ldots, n\}$ for all $n \in \mathbb{N}_0$.

Spectral identification allows to infer global information on the network structure from local measurements. It aims at estimating the spectrum $\sigma(L)$ of the Laplacian matrix $L$ from state measurements at a small number $p \ll n$ of nodes. More precisely, states measurements are used to compute the spectrum of the underlying Koopman operator through the DMD algorithm (see Schmid (2010); Tu et al. (2014)). Since the equilibrium $x^*$ is a stable equilibrium, the Koopman spectrum is directly related to the spectrum $\sigma(J)$ of the Jacobian matrix of the vector field (at the equilibrium)

$$J = I_n \otimes A - L \otimes BC^T \quad (2)$$

* This work is supported by the Namur Institute For Complex Systems (naXys) at University of Namur.

where $I_n$ is the identity matrix of order $n$, $A = \frac{\partial F}{\partial x}(x^*)$, $B = G(x^*)$ and $C = \frac{\partial H}{\partial x}(x^*)$. The spectral identification problem therefore boils down to estimating the Laplacian spectrum $\sigma(L)$ from the spectrum $\sigma(J)$ of the Jacobian matrix.

In Mauroy and Hendrickx (2017), the authors investigated the relationship between the spectrum of the Jacobian matrix $J$ and the spectrum of the Laplacian matrix $L$. It was shown that for all $\lambda \in \sigma(L)$, it exists an eigenvector $w$ such that
$$(A - \lambda BC^T)w = \mu w \qquad (3)$$
where $\mu \in \sigma(J)$ is an eigenvalue of $J$. In the specific case $r = 1$, they also proved that the relationship between the two spectra is one-to-one, so that one can retrieve the spectrum of $L$ from the spectrum of $J$. However the case $r = 1$ implies that the input and output signals $u_k$ and $y_k$ in (1) are scalar-valued, a condition that is not satisfied in many cases (see *e.g.* reaction-diffusion networks). The main contribution of this paper is to provide a generalized framework for spectral identification, which is valid for the case $r > 1$.

## 3. PROBLEM FEASIBILITY

In the following result, we show that the spectral identification problem is feasible under a mild assumption based on spectral moments. In particular, we define the $k$-th spectral moment of $L$ by
$$\mathcal{M}_k(L) = \frac{1}{n}\sum_{i=1}^{n} \lambda_i^k = \frac{1}{n}\operatorname{tr}(L^k) \qquad (4)$$
where $\operatorname{tr}(L)$ denotes the trace of $L$.

*Proposition 1.* If $\mathcal{M}_k(BC^T) \neq 0$ for all $k \in \mathcal{I}_n$, then the spectral identification problem is feasible, that is,
$$\sigma(J_1) = \sigma(J_2) \iff \sigma(L_1) = \sigma(L_2)$$
with $J_1 = I_n \otimes A - L_1 \otimes BC^T$ and $J_2 = I_n \otimes A - L_2 \otimes BC^T$.

**Proof.** Following similar lines as in Mauroy and Hendrickx (2017), we have, for $k \in \{1, \dots, n\}$
$$\mathcal{M}_k(J) = \frac{1}{mn} \operatorname{tr}\left[(I_n \otimes A - L \otimes BC^T)^k\right]$$
$$= \frac{1}{mn} \sum_{j=0}^{k} \binom{k}{j} \operatorname{tr}\left[(-L)^j \otimes A^{k-j} (BC^T)^j\right]$$
$$= (-1)^k \mathcal{M}_k(BC^T) \mathcal{M}_k(L)$$
$$+ \frac{1}{m}\sum_{j=0}^{k-1}\binom{k}{j}(-1)^j \operatorname{tr}\left[A^{k-j}(BC^T)^j\right]\mathcal{M}_j(L)$$
by applying the binomial theorem with $\binom{k}{j} = \frac{k!}{j!(k-j)!}$ and developing the traces of Kronecker-products.

This linear system of equations for the variables $M_k(L)$ is associated with a lower-triangular matrix that is invertible if its diagonal entries are all non zero, *i.e.* $(-1)^k \mathcal{M}_k(BC^T) \neq 0$ for all $k \in \{1, \dots, n\}$. In this case, the relationship between the first $n$ spectral moments $\mathcal{M}_k(L)$ and the first $n$ spectral moments $\mathcal{M}_k(J)$ is one-to-one. Since $\sigma(L)$ is uniquely determined by the $n$ first spectral moments $\mathcal{M}_k(L)$, it follows that $\sigma(J_1) = \sigma(J_2) \Rightarrow \sigma(L_1) = \sigma(L_2)$. The other implication directly follows from (3). ∎

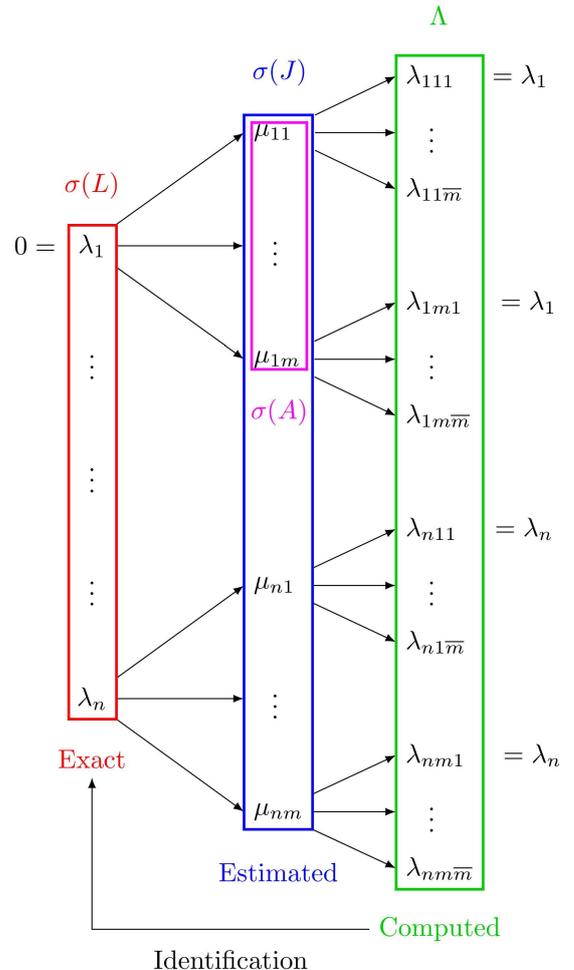

Fig. 1. Relationship between exact and identified eigenvalues.

If the assumptions of the above result are satisfied, the spectral moments $\mathcal{M}_k(L)$ can be obtained by solving the linear system of equations:
$$\mathcal{M}_k(L) = \frac{(-1)^k}{\mathcal{M}_k(BC^T)}\left(\mathcal{M}_k(J) + \right.$$
$$\left. \frac{1}{m}\sum_{j=0}^{k-1}\binom{k}{j}(-1)^{j+1}\mathcal{M}_j(L)\operatorname{tr}(A^{k-j}(BC^T)^j)\right). \qquad (5)$$

However, computing the $n$ eigenvalues of $L$ from the first $n$ spectral moments might lead to significant numerical error and should be avoided. In the next section, we propose an alternative method based on the properties of the characteristic polynomial
$$\chi(\lambda, \mu) = |A - \lambda BC^T - \mu I_m|$$
associated with the generalized eigenvalue problem (3).

## 4. SOLVING THE SPECTRAL NETWORK IDENTIFICATION PROBLEM

The relationship between the Laplacian spectrum $\sigma(L)$ and the spectrum of the Jacobian matrix $\sigma(J)$ can be described as follows. According to (3), any Laplacian eigenvalue $\lambda_i$, $i \in \mathcal{I}_n$, is in fact associated to $m$ eigenvalues of the Jacobian matrix $J$, which we denote by $\mu_{ij}, j \in \mathcal{I}_m$. In the context of the network identification problem, for each $\mu_{ij} \in \sigma(J)$ (estimated from data), one has to solve the generalized eigenvalue problem

$$(A - \mu_{ij} I_m) w = \lambda B C^T w \tag{6}$$

obtained by reorganizing (3). The solution can be computed by using the generalized Schur decomposition of $A$ and $BC^T$ (also called the QZ-algorithm, see Golub and Van Loan (1996)). However, there are $\overline{m}$ solutions $\lambda_{ijk}$, $k \in \mathcal{I}_{\overline{m}}$ for each value $\mu_{ij}$, where $\overline{m} \leq \min(r, m)$ is the rank of $BC^T$, one of which is $\lambda_i$. It follows that the set of computed eigenvalues $\Lambda = \{\lambda_{ijk} \mid i \in \mathcal{I}_n, j \in \mathcal{I}_m, k \in \mathcal{I}_{\overline{m}}\}$ is typically *larger* than $\sigma(L)$, i.e. $\sigma(L) \subsetneq \Lambda$, whenever $\overline{m} > 1$ (see Figure 1). The spectral identification problem therefore aims at identifying the Laplacian eigenvalues from those in $\Lambda$. Without loss of generality, we will denote these eigenvalues by $\lambda_{ij1} = \lambda_i$ for all $i \in \mathcal{I}_n$ and $j \in \mathcal{I}_m$.

One can also note that, for $\overline{m} = 1$, the solution of (6) is unique so that $\sigma(L) = \Lambda$. This particular situation is similar to the one in Mauroy and Hendrickx (2017), where $r = 1$ implies $\overline{m} = 1$.

### 4.1 Properties of the generalized eigenvalue problem

One can see that the exact computed eigenvalues $\lambda_{ij1} = \lambda_i$, with $j \in \mathcal{I}_m$, are repeated with multiplicity $m$. This property can be exploited unless spurious eigenvalues $\lambda_{ijk}$, with $k \neq 1$, are also repeated. In the worst case, the set of computed Laplacian eigenvalues may consist of $\overline{m}$ clusters of $n$ eigenvalues with multiplicity $m$, which we call indistinguishable clusters.

We will study the property of repeated eigenvalues by deriving two results on the characteristic polynomial $\chi(\lambda, \mu) = |A - \lambda BC^T - \mu I_m|$. Note that it can be rewritten as

$$\chi(\lambda, \mu) = \sum_{p=1}^{\overline{m}} \alpha_p \lambda^p - \sum_{q=0}^{m} \beta_q \mu^q + \sum_{s=1}^{\widetilde{m}} \sum_{t=1}^{m-s} \gamma_{st} \lambda^s \mu^t, \tag{7}$$

with $\widetilde{m} = \min(\overline{m}, m - 1)$.

**Lemma 2.** If $\text{tr}(BC^T) \neq 0$, then $\chi(\lambda, \mu)$ admits at least one term of the form $\lambda \mu^{m-1}$.

**Proof.** The characteristic polynomial $\chi(\lambda, \mu)$ is given by

$$\chi(\lambda, \mu) = |A - \lambda BC^T - \mu I_m|$$
$$= \sum_{\pi \in P_m} \varepsilon(\pi) \prod_{i=1}^{m} \left[ a_{i\pi(i)} - \lambda (BC^T)_{i\pi(i)} - \mu \delta_{i\pi(i)} \right], \tag{8}$$

where $P_m$ is the set of permutations of $m$ elements and $\varepsilon(\pi)$ is the signature of the permutation $\pi$.

Since $\delta_{i\pi(i)} = 0$ for all $i \neq \pi(i)$, terms of the form $\lambda^s \mu^t$, with $s + t \leq m$, are obtained only with permutations $\pi$ that admit at least $t$ fixed points. Furthermore, the only permutation in $P_m$ that has $m - 1$ fixed points is the identity, which is associated with the term

$$\prod_{i=1}^{m} \left[ a_{ii} - \lambda (BC^T)_{ii} - \mu \delta_{ii} \right],$$

in (8). By distributing, we can identify the terms in $\lambda \mu^{m-1}$, which add up and yield

$$(-1)^m \text{tr}(BC^T) \lambda \mu^{m-1}. \tag{9}$$

This term is nonzero by assumption. ∎

**Lemma 3.** Assume that $\overline{m} = m$ and that $\exists i \in \mathcal{I}_n$ such that, $\forall j_1 \neq j_2 \in \mathcal{I}_m$:

- $\mu_{ij_1} \neq \mu_{ij_2}$ are nonzero,
- $\{\lambda_{ij_1 k}\} = \{\lambda_{ij_2 k}\}$ (*i.e.* indistinguishable clusters).

Then the characteristic polynomial $\chi(\lambda, \mu)$ does not admit cross terms of the form $\mu^s \lambda^t$ with $s > 0$ and $t > 0$.

**Proof.** For fixed $i \in \mathcal{I}_n$ and $j \in \mathcal{I}_m$, we have that $\chi(\lambda_i, \mu_{ij}) = 0$ by definition of $\lambda_i$ and $\mu_{ij}$. It follows from (7) that

$$\sum_{p=1}^{\overline{m}} \alpha_p \lambda_i^p + \sum_{s=1}^{\widetilde{m}} \sum_{t=1}^{m-s} \gamma_{st} \lambda_i^s \mu_{ij}^t = \sum_{q=0}^{m} \beta_q \mu_{ij}^q.$$

For $k \in \mathcal{I}_{\overline{m}}$, $\lambda_{ijk}$ are the roots of

$$\chi(\lambda, \mu_{ij}) = \sum_{p=1}^{\overline{m}} \alpha_p \lambda^p + \sum_{s=1}^{\widetilde{m}} \sum_{t=1}^{m-s} \gamma_{st} \lambda^s \mu_{ij}^t - \sum_{q=0}^{m} \beta_q \mu_{ij}^q$$

$$= \sum_{p=1}^{\overline{m}} \alpha_p \left( \lambda^p - \lambda_i^p \right) + \sum_{s=1}^{\widetilde{m}} \sum_{t=1}^{m-s} \gamma_{st} \mu_{ij}^t \left( \lambda^s - \lambda_i^s \right). \tag{10}$$

Since $\overline{m} = m$, we have $\widetilde{m} = \min(\overline{m}, m - 1) = m - 1$ and (10) can be rewritten as

$$\chi(\lambda, \mu_{ij}) = \alpha_m \left( \lambda^m - \lambda_i^m \right) + \sum_{p=1}^{m-1} \left[ \left( \alpha_p + \sum_{t=1}^{m-p} \gamma_{pt} \mu_{ij}^t \right) \left( \lambda^p - \lambda_i^p \right) \right], \tag{11}$$

which is a polynomial in $\lambda$ that we denote by $\widetilde{\chi}_{ij}(\lambda)$. Furthermore, the assumption $\{\lambda_{ij_1 k}\} = \{\lambda_{ij_2 k}\}$ for all $j_1 \neq j_2 \in \mathcal{I}_m$ implies that the $m$ polynomials $\widetilde{\chi}_{ij}(\lambda)$ (*i.e.* for all $j \in \mathcal{I}_m$) share the same roots so that they are equal up to a multiplicative constant. Since they have the same coefficient $\alpha_m$, they are equal and the values $\sum_{t=1}^{m-p} \gamma_{pt} \mu_{ij}^t$ in (11) do not depend on $j$, that is, $\sum_{t=1}^{m-p} \gamma_{pt} \mu_{ij}^t = \xi_{pi}$.

Therefore, for each $p \in \mathcal{I}_{m-1}$, the coefficients $\gamma_{pt}$ are solutions of

$$\begin{pmatrix} \mu_{i1} & \mu_{i1}^2 & \cdots & \mu_{i1}^{m-p} \\ \vdots & \vdots & & \vdots \\ \mu_{im} & \mu_{im}^2 & \cdots & \mu_{im}^{m-p} \end{pmatrix} \begin{pmatrix} \gamma_{p1} \\ \vdots \\ \gamma_{p,m-p} \end{pmatrix} = \begin{pmatrix} \xi_{pi} \\ \vdots \\ \xi_{pi} \end{pmatrix}.$$

The matrix on the left-hand side is a Vandermonde matrix which is full rank since $\mu_{ij}$ are distinct and nonzero. It follows that the only solution is $\xi_{pi} = 0$ and $\gamma_{pt} = 0$, and therefore $\chi(\lambda, \mu)$ does not admit cross terms. ∎

The above result shows that, if the characteristic polynomials admit cross terms, then the computed eigenvalues will not consist in indistinguishable clusters of repeated eigenvalues. Conversely, if the polynomial does not admit cross terms, (10) writes $\chi(\lambda, \mu_{ij}) = \sum_{p=1}^{\overline{m}} \alpha_p \left( \lambda^p - \lambda_i^p \right)$ for all $i \in \mathcal{I}_n$ and $j \in \mathcal{I}_m$. It is clear that the polynomial does not depend on $j$, so that its roots $\{\lambda_{ijk}\}$ do not depend on $j$ either. In this case, computed eigenvalues consist of indistinguishable clusters of repeated eigenvalues.

## 4.2 Spectral network identification method

Our proposed method for spectral network identification is based on (i) the DMD algorithm to estimate the eigenvalues of the Jacobian matrix and (ii) the results of Section 4.1 to infer the Laplacian eigenvalues. It is summarized in Algorithm 1. Regarding the DMD algorithm, we follow the same guide lines as in Mauroy and Hendrickx (2017) (see Section 5 below). In particular, the method relies on the fact that each Laplacian eigenvalue $\lambda_i$ is obtained $m$ times in the computed set $\Lambda$ (one related to each value $\mu_{ij}$). Therefore, one can identify the Laplacian spectrum $\sigma(L)$ as the set of values in $\Lambda$ that are repeated with multiplicity $m$. Note that, due to numerical error in the estimation of the values $\mu_{ij}$, some tolerance $\varepsilon > 0$ should be considered to decide whether an eigenvalue is repeated or not, i.e. $\lambda$ is repeated if $|\lambda - \tilde{\lambda}| < \varepsilon$ for some $\tilde{\lambda} \in \Lambda$. Moreover, one has to ensure that spurious eigenvalues are not also repeated. This is verified in most cases under the necessary condition $\mathrm{tr}(BC^T) \neq 0$. Indeed, this condition implies that the characteristic polynomial admits cross terms (Lemma 2) and thus that the situation of $\overline{m}$ indistinguishable clusters of repeated eigenvalues is ruled out (Lemma 3), provided that $\overline{m} = m$ and the eigenvalues of $J$ are distinct and nonzero.

---

**Algorithm 1** Spectral network identification
---

**Inputs:** Dynamics $(A, B, C)$, state measurements (snapshot data), and tolerance $\varepsilon > 0$.
**Output:** Spectrum of the Laplacian matrix $\sigma(L)$.

**Eigenvalue estimation:**
1: Estimate the eigenvalues $\mu \in \sigma(J)$ from state measurements using the DMD algorithm.

**Computing:**
2: **for** each $\mu \in \sigma(J)$ **do**
3:     Solve (6) using the QZ-algorithm and store the solutions in $\Lambda$.
4: **end for**

**Filtering:**
5: Count the multiplicity of each value in $\Lambda$ with the tolerance $\varepsilon$, without pruning the values already considered for the multiplicity of other values.

**Averaging:**
6: **for** each value of $\Lambda$ with multiplicity $m$ **do**
7:     Define $\lambda_i$ as the average of the $m$ repeated eigenvalues and store $\lambda_i$ in $\sigma(L)$.
8: **end for**

---

## 4.3 Illustration with random networks

We consider a random network of $n = 50$ nodes with an adjacency matrix $W$ whose entries are uniformly distributed in $[0, 1]$. The linearized dynamics is characterized by matrices $A$ and $B$ with entries that are uniformly distributed in $[0, 20]$ (with $m = 2$) and by a matrix $C = I_2$. The spectrum of the Jacobian matrix is used as input for Algorithm 1 with $\varepsilon = 10^{-4}$. In this case, spectral identification is feasible since $\overline{m} = m = 2$ and $\mathrm{tr}(BC^T) \neq 0$. Indeed,

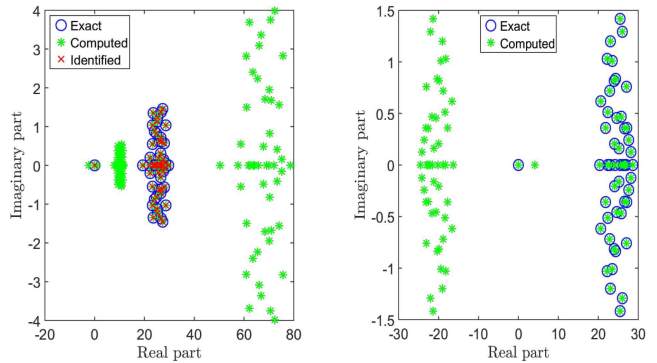

(a) Exact spectral identification with $\mathrm{tr}(BC^T) \neq 0$.

(b) Indistinguishable clusters with $\mathrm{tr}(BC^T) = 0$.

Fig. 2. Blue circles represent the exact Laplacian spectrum $\sigma(L)$, green stars represent the computed eigenvalues $\Lambda$ and red crosses represent the identified Laplacian spectrum.

we observe only one cluster of repeated eigenvalues, corresponding to correct eigenvalues (Figure 2(a)). In contrast, in the case where $\mathrm{tr}(BC^T) = 0$ (imposing $B_{11} = 1$ and $B_{22} = -1$), we observe two (indistinguishable) clusters of repeated eigenvalues (Figure 2(b)). Note that in the latter case, one can still identify correct eigenvalues by exploiting the fact that Laplacian eigenvalues are positive.

## 5. NUMERICAL EXPERIMENTS

In this section, our spectral network identification method is illustrated with two examples.

### 5.1 Data and DMD algorithm

We consider Erdős-Rényi networks of size $n$, where each entry of the adjacency matrix $W$ has a probability of $p_{\mathrm{edge}}$ to be uniformly distributed in $[0, 1]$ and $1 - p_{\mathrm{edge}}$ to be 0. Note that $n$ and $p_{\mathrm{edge}}$ are supposedly unknown.

We consider both linear and nonlinear dynamics:

(1) Linear dynamics
$$\begin{cases} \dot{x}_k = \begin{pmatrix} -1 & -2 \\ 1 & -1 \end{pmatrix} x_k + \begin{pmatrix} 1 & 0 \\ 0 & 2 \end{pmatrix} u_k \\ u_k = \sum_{j=1}^n W_{kj}(y_k - y_j) \\ y_k = x_k. \end{cases} \quad (12)$$

(2) Nonlinear Brusselator dynamics (Prigogine (1980))
$$\begin{cases} \dot{x}_1 = 1 - (b+1)\, x_1 + a\, x_1^2\, x_2 \\ \dot{x}_2 = b\, x_1 - a\, x_1^2\, x_2, \end{cases} \quad (13)$$

where $a, b \geq 0$ are parameters. We set $a = 15$ and $b = 9$ so that the fixed point is $x^* = (1, b/a) = (1, 0.6)$. The coupling is defined by $B = \mathrm{diag}(1, 4.5)$ and $C = I_2$.

We generate $q$ trajectories $X^{(j)}$ $(j = 1, \ldots, q)$ from different initial conditions whose components are uniformly distributed in an interval of length $10^{-4}$ and centered on the associated component of $x^*$. We take $K+1$ snapshots data on each trajectory with a sampling period $\Delta t$, $X^{(j)}(k\Delta t)$ ($k = 0, \ldots, K$). We consider the measurement function $f : \mathbb{R}^{mn} \to \mathbb{R}$ which returns only the first state of the first node. The DMD algorithm (see Schmid (2010); Tu et al. (2014)) is used with the data measurements $f(X^{(j)}(k\Delta t))$ and $c - 1$ delayed sequences $f(X^{(j)}((k+1)\Delta t))$, ..., $f(X^{(j)}((k+c-1)\Delta t))$ (see Mauroy and Hendrickx (2017) for more details). The parameters used in each experiment are summarized in Table 1.

| **Case** | $n$ | $p_{\text{edge}}$ | $q$ | $K$ | $\Delta t$ | $c$ | $\delta$ | $\varepsilon$ |
|---|---|---|---|---|---|---|---|---|
| L | 10 | 0.65 | 10 | 25 | 0.4 | 2 | 5 | 0.05 |
| NL | 10 | 0.65 | 10 | 25 | 0.4 | 2 | 5 | 0.1 |
| NL | 100 | 0.3 | 10 | 50 | 0.4 | 2 | 5 | 1 |

Table 1. Parameters used for numerical experiments in the linear (L) and nonlinear (NL) cases.

### 5.2 Numerical results

In the linear case and with a small network ($n = 10$), Algorithm 1 is efficient to recover all Laplacian eigenvalues (Figure 3).

With nonlinear dynamics, the DMD algorithm is not able to compute accurately all the eigenvalues of the Jacobian matrix so that it becomes more difficult to capture repeated values in $\Lambda$. Similarly, it is not possible to capture all the eigenvalues of the Jacobian matrix in the case of large networks. For these cases, our method must be adapted.

In the case of small graphs, we adapt the filtering part in Algorithm 1 as follows. Assuming that the different (correct and spurious) clusters of computed eigenvalues in $\Lambda$ lie in separated regions of the complex plane, we can use repeated eigenvalues identified by Algorithm 1 to define a region of the complex plane where *all* eigenvalues—even not repeated—supposedly belong to the Laplacian spectrum. We will define this region as the convex hull of the eigenvalues identified by Algorithm 1, and will select all computed eigenvalues that lie inside this set with the tolerance $\varepsilon$.

Note that if only real eigenvalues are identified, then the convex hull is not defined. In this case, the approximation of the Laplacian spectrum is given by all $\lambda \in \Lambda$ such that $\text{Re}\{\lambda\} \in [0, \rho + \varepsilon]$ and $\text{Im}\{\lambda\} \in [-\varepsilon, \varepsilon]$, where $\rho$ is the identified eigenvalue with the largest real part. We use this adapted method to approximate the Laplacian eigenvalues in the case of the Brusselator dynamics on the same network as in the linear case (see Figure 4).

In the case of large networks, one can make the acceptable approximation that the Laplacian eigenvalues are uniformly distributed in a region of the complex plane. Again, this region is defined as the convex hull of repeated eigenvalues obtained with Algorithm 1. Then, the first two spectral moment of $L$ can be approximated using the area, centroid and second moment of area of the convex hull (see Mauroy and Hendrickx (2017) for more details). This method is used for a large network ($n = 100$) with the Brusselator dynamics (Figure 5) and the results are given

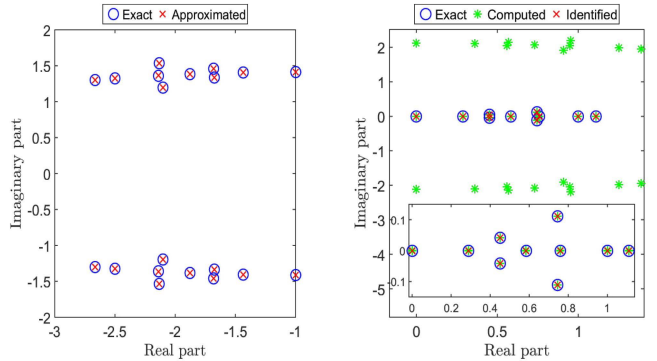

(a) Estimated eigenvalues of the Jacobian matrix.

(b) Computed Laplacian eigenvalues.

Fig. 3. The spectral network identification method recovers the Laplacian eigenvalues in the case of a network of $n = 10$ nodes, with the linear dynamics (12). Inset in panel (b): close-up view on the exact Laplacian eigenvalues.

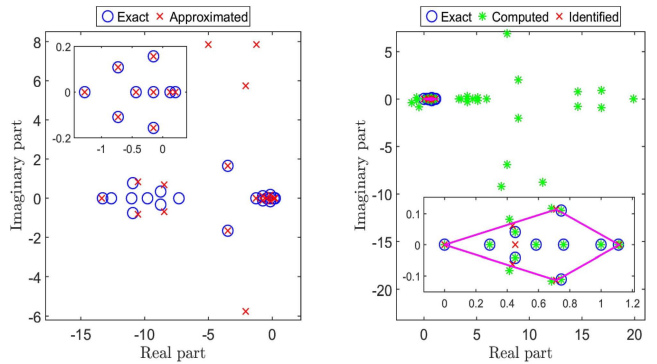

(a) Estimated eigenvalues of the Jacobian matrix.

(b) Computed Laplacian eigenvalues.

Fig. 4. The adapted spectral network identification method is used in the case of a network of size $n = 10$ with the Brusselator dynamics (13). The magenta line is the convex hull. Inset in panel (b): close-up view on the exact Laplacian eigenvalues.

| $\boldsymbol{n = 100}$ | Estimated | Exact |
|---|---|---|
| $\mathcal{M}_1(L)$ | 5.55 | 5.05 |
| $\mathcal{M}_2(L)$ | 35.04 | 29.48 |
| $\lambda_2$ | 0.29 | 1.45 |
| $\lambda_n$ | 8.75 | 9.85 |

Table 2. Spectral identification of a network of size $n = 100$ with the Brusselator dynamic.

in Table 2. Note that $\lambda_2$ and $\lambda_n$ are taken as the minimal and maximal real part of the convex hull, respectively.

### 6. CONCLUSION

This paper considers the spectral network identification framework introduced in Mauroy and Hendrickx (2017) and extends its applicability to the case of vector-valued input and output coupling signals. The feasibility problem has been investigated, which turns out to be more involved

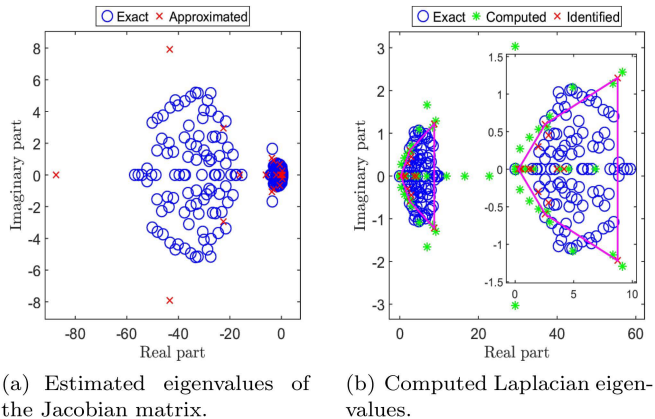

(a) Estimated eigenvalues of the Jacobian matrix.

(b) Computed Laplacian eigenvalues.

Fig. 5. The adapted spectral network identification method is used in the case of a large network of size $n = 100$ with the Brusselator dynamic (13). The magenta line is the convex hull. Inset in panel (b): close-up view on the exact Laplacian eigenvalues.

in this case. Moreover we have proposed a method for solving the spectral identification problem, which is based on the properties of a generalized eigenvalue problem. This method has been illustrated with several numerical experiments, which showed that all Laplacian eigenvalues can be inferred in the case of linear dynamics and small networks, or spectral moments and bounds on Laplacian eigenvalues can be obtained otherwise.

This work leads to several research perspectives. First of all, the proposed method is not robust to numerical errors on the eigenvalues estimated by the DMD algorithm. We have proposed a heuristic solution to this issue, which could be improved. In this context, specific extensions of the DMD algorithm could also be used (*e.g.*, Brunton et al. (2013), Williams et al. (2015), Hemati et al. (2017), Li et al. (2017), Gulina and Mauroy (2021)). Moreover, the results and methods developed in this work are based on the assumption of a network of identical units. Further research work could investigate the case of nonidentical units with possibly several equilibria or even admitting a limit-cycle. Finally, additional work could be needed to infer relevant network properties from identified features of the Laplacian spectrum. While the spectra of the adjacency and modularity matrices of random graphs have been studied in Nadakuditi and Newman (2013), the reverse problem of estimating distributions of expected degree from the identified Laplacian spectrum could be addressed.

## ACKNOWLEDGEMENTS

The authors are grateful to Germain van Bever and Cédric Simal for fruitful discussions.